\documentclass[preprint,12pt]{elsarticle}
\usepackage{algorithm}
\usepackage{algpseudocode}
\usepackage{esvect}
\usepackage[title]{appendix}
\usepackage{mathrsfs}
\usepackage{array,siunitx}
\usepackage{amsfonts, amsmath}
\usepackage{amsmath}
\usepackage{graphicx}
\usepackage{color}

\usepackage{bm}
\usepackage{caption}
\usepackage{verbatim}
\usepackage{hyperref}
\usepackage{bmpsize}
\usepackage{float}
\restylefloat{table}
\usepackage[table, dvipsnames]{xcolor}
\usepackage[margin=0.9in]{geometry}
\usepackage[belowskip=-15pt,aboveskip=2pt]{caption}
\usepackage{makecell}
\usepackage{cite}
\usepackage{relsize}
\allowdisplaybreaks

\newtheorem{theorem}{Theorem}
\newtheorem{lemma}{Lemma}
\newtheorem{remark}{Remark}
\newproof{pf}{Proof}
\newcommand{\be}{\begin{equation}}
	\newcommand{\ee}{\end{equation}}
\newcommand{\bea}{\begin{eqnarray}}
	\newcommand{\eea}{\end{eqnarray}}
\newcommand{\beas}{\begin{eqnarray*}}
	\newcommand{\eeas}{\end{eqnarray*}}

\def\qed{\hbox{\vrule width 6pt height 6pt depth 0pt}}
\newcommand{\vertiii}[1]{{\left\vert\kern-0.25ex\left\vert\kern-0.25ex\left\vert #1
		\right\vert\kern-0.25ex\right\vert\kern-0.25ex\right\vert}}

\newcommand{\normiii}[1]{{\left\vert\kern-0.25ex\left\vert\kern-0.25ex\left\vert #1
		\right\vert\kern-0.25ex\right\vert\kern-0.25ex\right\vert}}

\usepackage{amssymb}

\journal{Journal of Mathematical Analysis and Applications}

\begin{document}

\begin{frontmatter}



\title{Uniqueness of Steady Navier-Stokes under Large Data by Continuous Data Assimilation}


\author{Xuejian Li}

\affiliation{organization={School of Mathematics and Statistics Sciences},
            addressline={Martin Hall}, 
            city={Clemson},
            state={SC},
            postcode={29634},
       country={USA} } 
            

\begin{abstract}
We propose a continuous data assimilation (CDA) method to address the uniqueness 
problem for steady Navier-Stokes equations (NSE). The CDA method incorporates spatial observations into the NSE, and we prove that with sufficient observations, the CDA-NSE system is well-posed even for large data where multiple solutions may exist.  
This CDA idea is in general helpful to determine solution for non-uniqueness partial differential equations (PDEs).

\end{abstract}




\end{frontmatter}


\section{Introduction}
The Navier-Stokes equations (NSE) are fundamental in modeling fluid mechanics. 
On $\mathbb{R}^d, d=2,3$, the steady NSE for incompressible Newtonian fluids is given by
\begin{equation}\label{NS}
	\left\{\begin{aligned}
		&-\nu \Delta u+(u\cdot\nabla) u+ \nabla p={f} \quad \text{in}~\Omega,\\
		&\nabla\cdot {u}=0\quad \text{in}~\Omega,\\
		&{u}=0 \quad \text{on}~~\partial\Omega,
	\end{aligned}\right.
\end{equation}
where 
$u$ is the velocity of fluid, $p$ is the kinetic pressure, $\nabla\cdot u=0$ indicates that the fluid is incompressible, 
${f}$ is the external force, and $\nu$ is the viscosity of the fluid. The parameter $Re=\frac{1}{\nu}$ plays the role of Reynolds number.

It is well-known that for small data, i.e. small $Re$ and $f$, there exists a unique solution for the system (\ref{NS}). However, while $Re$ or $f$ increases and crosses certain critical bounds, the NSE can lose uniqueness and admits multiple solutions that fall into different branches \citep{MR548867}. This phenomenon is often encountered in practice, and these non-unique solutions are often called isolated solutions or branches of nonsingular solutions \citep{MR548867,MR2442411}. Numerically finding such solutions is especially difficult due to non-uniqueness making nonlinear iterative solvers less effective.

The main interest of this paper is showing that using continuous data assimilation (CDA)  \citep{MR3994475,MR2036872,MR3183055,MR3319381,MR4264472} can overcome the uniqueness difficulty for the steady NSE.  While CDA is generally used with time dependent problems, the type of nudging employed by CDA can also be applied to steady problems, however the notion of continuous (in time) is no longer valid; still, we refer to it as CDA in this paper.  To define the steady CDA-NSE system, let $I_Hu$ represent an interpolant operator (or observation operator) based on spatial observations of a NSE solution $u$ of system (\ref{NS})
at a coarse resolution mesh size $H$ (requirements for $I_H$ are given in section 2). To uniquely identify the solution for system (\ref{NS}) associated with the measurements $I_Hu$, we propose the following CDA-NSE system:
\begin{equation}\label{NSCDA}
	\left\{\begin{aligned}
		&-\nu \Delta w+(w\cdot\nabla) w+ \nabla z+\mu(I_Hw-I_Hu)={f} \quad \text{in}~\Omega,\\
		&\nabla\cdot {w}=0\quad \text{in}~\Omega,\\
		&{w}=0 \quad \text{on}~~\partial\Omega,
	\end{aligned}\right.
\end{equation}
where $\mu(I_Hw-I_H u)$ is a nudging term driving state $w$ towards to the observations, and $\mu$ is a positive relaxation parameter that emphasizes the observations accuracy. In this  context, we consider accurate spatial observations, and thus there are no size restrictions on $\mu$.

We show that with enough observations, i.e. that $H$ is sufficiently small, the CDA-NSE (\ref{NSCDA}) has a unique solution even for large $Re$ and $f$, and the CDA-NSE solution is identical to the isolated NSE solution that corresponds to the observed state.   The analysis and results in this paper may have a positive influence on developing effective iterative solvers for the steady NSE with large Reynolds number or external forces when observations are available.

While this note studies the NSE, a similar idea can lead to wellposedness for related steady multi-physics problems, such as magnetohydrodynamics or Boussinesq systems.

\section{Uniqueness analysis}
Before formally presenting the main results, we briefly introduce necessary preliminaries.  Consider $\Omega$ as an open bounded domain, denote the natural function spaces by 
\begin{align}
	&Q:=\{v\in {L}^2(\Omega): \int_{\Omega}vdx=0\},\\
	& X:=
	\{v\in H^1\left(\Omega\right): v=0~~\text{on}~ \partial\Omega\}, \\
	& V:=
	\{v\in X:  \nabla \cdot v=0\}.
\end{align}
Let $(\cdot,\cdot)$ denote the $L^2(\Omega)$ inner product that induces the $L^2$ norm $\|\cdot\|$,  $H^{-1}$ and $V^*$ denote the dual spaces of $X$ and $V$, respectively. In addition, let $\langle z,v\rangle_{-1}$ denote the action of $z\in H^{-1}$ on $v\in X$ and  $\langle z,v\rangle_{*}$ denote the action of $z\in V^*$ on $v\in V$, respectively. Also, 
\begin{align*}
	\|z\|_{-1}=\sup_{\forall v\in X}\frac{\langle z,v\rangle_{-1}}{\|\nabla v\|},~~\|z\|_{*}=\sup_{\forall v\in V}\frac{\langle z,v\rangle_{*}}{\|\nabla v\|}.
\end{align*}
 The weak form of NSE (\ref{NS}) is to find $(u,p)\in X\times Q$ such that 
\begin{equation}
	\label{weakNavier}
	a\left({u},v\right)+b\left({u},{u},v\right)+(p,\nabla\cdot v)=\left\langle{f},v\right\rangle_{-1}~~\forall v\in X,~~(\nabla\cdot{u},q)=0~~\forall q\in Q,
\end{equation}
where $a(\cdot,\cdot)$ and $b(\cdot,\cdot,\cdot)$ are defined as follows:
\begin{align*}
	&	a(u,v)=\left(\nu\nabla u, \nabla v\right)~~\forall u,v\in X\\&b(u,w,v)=\left((u\cdot\nabla)w, v\right)~~\forall u,w,v\in X.
\end{align*}
Note that due to inf-sup condition holding on $X\times Q$ \citep{MR548867,constantinfoiasciprian1989}:
\begin{align*}
	\inf\limits_{0\neq q\in Q}\sup\limits_{ 0 \neq{v}\in X}\frac{(q,\nabla\cdot v)}{\left\|q\right\|_{Q}\left\|{v}\right\|_{{X}}}\geq\beta>0,
\end{align*}
the system (\ref{weakNavier}) is equivalent to: Find $u\in V $ satisfying
\begin{equation}\label{wd}
	a\left({u},v\right)+b\left({u},{u},v\right)=\left\langle{f},v\right\rangle_* ~~\forall v\in V.
\end{equation}
For the trilinear term $b(\cdot,\cdot,\cdot)$, the following inequalities hold \citep{MR2442411,MR603444}:
\begin{align}
		\label{3dd}	b(u,w,v)&\leq M\|\nabla u\|\|\nabla w\|\|\nabla v\|~~\text{for $d=2$ and $d=3$},\\
	\label{3d}	b(u,w,v)&\leq M_1\|u\|^{\frac{1}{2}}\|\nabla u\|^{\frac{1}{2}}\|\nabla w\|\|\nabla v\|~~\text{for $d=2$ and $d=3$},\\
	\label{2d}		b(u,w,v)&\leq M_2\|u\|^{\frac{1}{2}}\|\nabla u\|^{\frac{1}{2}}\|\nabla w\|\|v\|^{\frac{1}{2}}\|\nabla v\|^{\frac{1}{2}}~~\text{for $d=2$}.
\end{align}
Here, $M$, $M_1$, and $M_2$ are positive constants depending on $\Omega$.

We recall the classical  well-posedness result for equation (\ref{wd}) \citep{MR2442411,MR603444}:

\begin{lemma}
	Let $\alpha=M\nu^{-2}\|f\|_{*}$. For any $f\in V^*$ and $\nu$,	there exists at least one solution for NSE (\ref{wd}). Besides this, every solution of (\ref{wd}) satisfy a priori estimate 
	\begin{align}\label{Pri}
		\|\nabla u\|\leq \nu^{-1}\|f\|_{*}.
	\end{align}
	Furthermore, if $\alpha<1$, 
	the solution is unique.\\ 
\end{lemma}
The restriction $\alpha< 1$ is usually referred as the small data condition for steady NSE. In this same spirit, we refer  to $\alpha\geq 1$ as the case of large data. 

\vspace{1mm}
Given interpolated observations $I_Hu$,  
the weak form of the CDA-NSE (\ref{NSCDA}) is to find $w\in V$ such that 
\begin{equation}\label{cdaweak}
	a\left({w},v\right)+b\left({w},{w},v\right)+\mu (I_Hw-I_H  u,I_Hv)=\left\langle{f},v\right\rangle_* ~~\forall v\in V. 
\end{equation}
\begin{remark}
Note that in (\ref{cdaweak}), $\mu (I_H  u,v)=\mu (I_H  u,I_Hv)~\forall u,v\in X$ in the case that $I_H$ is the $L^2$ projection onto the coarse mesh space, and for general $I_H$ that all results below still hold if you used $\mu(I_H u,v)$ instead of $\mu (I_H  u,I_Hv)$ but there would be stronger restrictions on $\mu$ and $H$.
\end{remark}
 In the remainder of the paper, we assume the interpolant $I_H$ is linear and have the properties:
\begin{align}
	&\label{interpolationi}	\|I_Hv-v\|\leq C_IH\|\nabla v\|,~~\|I_Hv\|\leq C\|v\|~~\forall v\in X.
\end{align}
Such interpolant generally exists in finite approximation theory, for instance the $P_1$ finite element interpolation\citep{MR4264472}:
\begin{align*}
		I_Hv:=\sum_{j=1}^{N_1} v(x_{H}^{j})\phi_j~~\forall v\in X.
\end{align*}
Here, $H$ can be the finite element mesh size, $N_1$ is the number of finite element nodes, $x_{H}^{j}$ is the $j^{th}$ finite element node, and $\{\phi_j\}_{j=1}^{N_1}$ are the degree one polynomial finite element basis. 
 
Based on Leray-Schauder fixed point theorem\footnote{This is the only place where the inequality $\|I_Hv\|\leq C\|v\|$ in (\ref{interpolationi}) is in need.}, it is not difficult to prove the CDA-NSE (\ref{cdaweak}) has at least one solution for any non-negative $\mu$ and $H$.  Additionally, one can observe if $w=u$ is a solution to (\ref{cdaweak}), then the existence is established this way as well. In the following, we focus the relation between equations (\ref{cdaweak}) and (\ref{wd}) and show the uniqueness of (\ref{cdaweak}). 

\begin{theorem}\label{equi}
	Assume $f\in V^*$ and $u$ is a solution of (\ref{wd}). If $\alpha<1$, for any given $H$ and $\mu$, the CDA-NSE (\ref{cdaweak}) is equivalent to the NSE (\ref{wd}) in sense that the solution $w$ to (\ref{cdaweak}) is unique and equal to $u$. If $\alpha\geq 1$, under the  condition 
	\begin{equation}\label{equic}
		\begin{split}
		H\leq \frac{2M^2}{3\sqrt{3}C_IM_1^2\alpha^2}~~\text{and~ } \mu\geq \frac{\nu}{4C_I^2H^2},
		\end{split}
	\end{equation}
the CDA-NSE (\ref{cdaweak}) has a unique solution which is exactly the isolated solution of NSE (\ref{wd}) that corresponds to the observed state, that is, we also have $w=u$.  
\end{theorem} 
\begin{pf}

	Subtracting equation (\ref{cdaweak}) from (\ref{wd}), we have
	\begin{equation}\label{M1}
		\begin{split}
			0&=a\left({w},v\right)-a\left(u,v\right)+b\left({w},{w},v\right)-b\left(u,u,v\right)+\mu (I_Hw-I_H u,I_Hv)\\
			&=a\left({w}-u,v\right)+b(w,w-u,v)+b(w-u,u,v)+\mu (I_Hw-I_H u,I_Hv).
		\end{split}
	\end{equation}
	Taking $v=w-u$, and using (\ref{3dd}) and (\ref{Pri}),  we obtain
	\begin{equation}\label{equ1}
		\begin{split}
			&\nu\|\nabla({w}-u)\|^2+\mu\|I_Hw-I_H u\|^2=-b(w-u,u,w-u)\\
			&\leq M\|\nabla(w-u)\|^2\|\nabla u\|\leq M\nu^{-1}\|f\|_{*}\|\nabla(w-u)\|^2.
		\end{split}
	\end{equation}
	Rearranging (\ref{equ1}) gives us 
	\begin{equation}
		\begin{split}
			\nu(1-M\nu^{-2}\|f\|_{*})\|\nabla({w}-u)\|^2+\mu\|I_Hw-I_H u\|^2&\leq 0.
		\end{split}
	\end{equation}
	If $\alpha<1$, it is clear to see $\|\nabla({w}-u)\|=0$ is always true, i.e., $w=u$. Thus with $\alpha<1$, the NSE (\ref{wd}) has a unique solution, and so $w=u$ is the unique CDA-NSE solution.

	Next, we consider the case $\alpha\geq1$. Continuing from the equality in (\ref{equ1}), using inequalities (\ref{3d}) and (\ref{Pri}) and generalized Young's inequality, we have
	\begin{equation}\label{equ2}
		\begin{split}
			&\nu\|\nabla({w}-u)\|^2+\mu\|I_Hw-I_H u\|^2=-b(w-u,u,w-u)\\
			&\leq M_1\|w-u\|^{\frac{1}{2}}\|\nabla(w-u)\|^{\frac{1}{2}}\|\nabla u\|\|\nabla(w-u)\|\\
			&\leq M_1\nu^{-1}\|f\|_{*}\|\nabla(w-u)\|^{\frac{3}{2}}\|w-u\|^{\frac{1}{2}}\\
			&   \le \frac{M_1}{M}\nu\alpha \|\nabla(w-u)\|^{\frac{3}{2}}\|w-u\|^{\frac{1}{2}} \\ 
			&  \le \frac{\nu}{2}\|\nabla(w-u)\|^2 + \frac{27M_1^4\nu \alpha^4}{32M^4} \|w-u\|^{2} .			
		\end{split}
	\end{equation}
	Applying inequality (\ref{interpolationi}) and the norm inequality $\frac{\|a-b\|^2}{2}\leq \|a-c\|^2+\|c-b\|^2$, we bound the left side of (\ref{equ2}) from below as
	\begin{equation}\label{equ3}
		\begin{split}
			&\nu\|\nabla({w}-u)\|^2+\mu\|I_Hw-I_H u\|^2\\
			&\geq \frac{3\nu}{4}\|\nabla({w}-u)\|^2+\frac{\nu}{4C_I^2H^2}\|({w}-u)-I_H({w}-u)\|^2+\mu\|I_Hw-I_H u\|^2\\
			&\geq\frac{3\nu}{4}\|\nabla({w}-u)\|^2+\frac{\lambda}{2}\|w- u\|^2,
		\end{split}
	\end{equation}
	where $\lambda=\min \{\frac{\nu}{4C_I^2H^2},\mu\}$. \\
	Combining (\ref{equ2}) and (\ref{equ3}) leads to 
	\begin{equation}\label{equ4}
		\begin{split}
			\frac{\nu}{4}\|\nabla({w}-u)\|^2+\left(\frac{\lambda}{2}-\frac{27M_1^4\nu \alpha^4}{32M^4}\right)\|w- u\|^2\leq 0.
		\end{split}
	\end{equation}
	Recall that $\mu$ can be large and there is no upper bound on $\mu$ that arises in our analysis, we thus consider $\mu$ large enough so that $\lambda= \frac{\nu}{4C_I^2H^2}$. If $\frac{\lambda}{2}-\frac{27M_1^4\nu \alpha^4}{32M^4}\geq 0$ is satisfied, i.e.,
	\begin{equation}\label{equ5}
		\begin{split}
			H\leq \frac{2M^2}{3\sqrt{3}C_IM_1^2\alpha^2},
		\end{split}
	\end{equation}
	then $\|\nabla({w}-u)\|=0$ holds. 	

	Finally, since the solutions of the steady NSE are isolated, then $w$ must be the observed isolated solution of equation (\ref{wd}) and thus is unique to (\ref{cdaweak}) as well. 
	This completes the proof. \qed
\end{pf}
\begin{remark}
	The condition on $H$ in (\ref{equic}) is less restrictive for $d=2$. Continuing from (\ref{equ1}), using inequality (\ref{2d}), (\ref{Pri}), and Young's inequality, we have
	\begin{equation}\label{equ22}
		\begin{split}
			&\nu\|\nabla({w}-u)\|^2+\mu\|I_Hw-I_H u\|^2=-b(w-u,u,w-u)\\
			&\leq M_2\|w-u\|^{\frac{1}{2}}\|\nabla(w-u)\|^{\frac{1}{2}}\|\nabla u\|\|\nabla(w-u)\|^{\frac{1}{2}}\|w-u\|^{\frac{1}{2}}\\
			&\leq \frac{M_2}{M}\nu\alpha\|\nabla(w-u)\|\|w-u\|\\&\leq \frac{\nu}{2}\|\nabla(w-u)\|^2+\frac{M_2^2\nu\alpha^2}{2M^2}\|w-u\|^2.
		\end{split}
	\end{equation}
	Combining (\ref{equ22}) and (\ref{equ3}) leads to 
	\begin{equation}\label{equ44}
		\begin{split}
			\frac{\nu}{4}\|\nabla({w}-u)\|^2+\left(\frac{\lambda}{2}-\frac{M_2^2\nu\alpha^2}{2M^2}\right)\|w- u\|^2\leq 0.
		\end{split}
	\end{equation}
	Consider $\lambda=\frac{\nu}{4C_I^2H^2}$. If $\frac{\lambda}{2}-\frac{M_2^2\nu\alpha^2}{2M^2}=\frac{\nu}{8C_I^2H^2}-\frac{M_2^2\nu\alpha^2}{2M^2}\geq 0$ is satisfied, i.e.,
	\begin{equation}\label{equ55}
		\begin{split}
			H\leq \frac{M}{2C_IM_2\alpha},
		\end{split}
	\end{equation}
	then $\|\nabla({w}-u)\|=0$ must hold . Similarly, condition (\ref{equ55}) is also sufficient for the uniqueness. Note that, compared to inequality (\ref{equic}), this is a significantly less restriction on $H$.
\end{remark}

\section{Conclusion}
We proposed a CDA-NSE alteration of the steady NSE system that incorporates observables through the CDA nudging process, and proved that with enough observables the system is well-posed for any data.  We showed a sufficient condition for how much observables is needed for well-posedness, and the amount scales with the size of the data.  The analysis and results in this paper provides a mathematical foundation for incorporating CDA into iterative nonlinear solvers for the steady NSE, which is a subject of ongoing research by the author.

\section*{Acknowledgments}
This work is partially supported by NSF Grant DMS 2152623.

 \bibliographystyle{elsarticle-num} 
\bibliography{Xuejian_Li_ref}
\end{document}